\documentclass[12pt]{article}
\usepackage{times}
\usepackage{amsmath}
\usepackage{amssymb}

\newtheorem{theorem}{Theorem}[section]
\newtheorem{proposition}[theorem]{Proposition}
\newtheorem{lemma}[theorem]{Lemma}

\newtheorem{definition}[theorem]{Definition}
\newtheorem{remark}[theorem]{Remark}
\def\SL{\mathop{\operatorname{\rm SL}}\nolimits}
\def\Hom{\mathop{\operatorname{\rm Hom}}\nolimits}  
\def\ind{\mathop{\operatorname{\rm Ind}}\nolimits}  
  
\def\aff{\mathop{\operatorname{\rm aff}}\nolimits}  
\def\End{\mathop{\operatorname{\rm End}}\nolimits}  
\def\Sym{\mathop{\operatorname{\rm Sym}}\nolimits}  
\def\infl{\mathop{\operatorname{\rm infl}}\nolimits}  
\begin{document}
\title{Quantized Algebras of Functions on\\ Affine Hecke Algebras\footnote{The work was supported in part by  National Foundation for Research in Fundamental Sciences, Vietnam, Alexander von Humboldt Foundation, Germany, and was completed during the visit of the author at the Department of Mathematics, The University of Iowa, U.S.A.}
}
\author{Do Ngoc Diep}
\maketitle
\begin{abstract}The so called quantized algebras of functions on affine Hecke algebras of type A and the corresponding $q$-Schur algebras are defined and their irreducible unitarizable representations are classified.
\end{abstract}
\section*{Introduction}
The algebras of functions on groups define the structure of the groups themselves: the algebras of continuous functions on topological groups define the structure of the topological groups. This essentially is the so called Pontryagin duality for Abelian locally compact groups and the Tannaka-Krein duality theory for compact groups.  The smooth functions on Lie groups define the structure of Lie groups. It is the essential fact that in this case we can produce the harmonique analysis on genral Lie groups.  The quantized algebras of functions on quantum groups defined the structure of quantum groups etc. In the same sense we define quantized algebras of functions which define the structure of quantum affine Hecke algebras. 
Let us discuss a little bites in more detail. Let us denote by $\mathfrak g$ a Lie algebra over the field of complex numbers, $U(\mathfrak g)$ its universal eveloping algebra, $\lambda \in P^*$ a positive highest weight, $V_v(\lambda)$ the associated representation of type I, i.e. with a positive defined Hermite form $(.,.)$ and $(x v_1, v_2)= (v_1. x^*v_2), \forall v_1, v_2 \in V_v(\lambda)$, of the  quantized universal enveloping algebra $U_v(\mathfrak g)$. Let $\{ v_\mu^\nu \}$ be an othogonal basis of $V_v(\lambda)$. 
Consider the matrix elements of the representation defined by 
$$C^\lambda_{\nu, s; \mu, r}(x) : = (xv_\mu^r, v_\nu^s)$$ and the linear span $\mathcal F_v(G) := \langle C^\lambda_{\nu,s; \mu, r} \rangle$. It was shown in L. Korogodski and Y. Soibelman \cite{korosoibel} that indeed it is equipped with a structure of an Hopf algebra, the so called the quantized algebra of functions on the quantum group corresponding to $G$. It was shown also that this algebra is generalized by the matrix coefficients of the standard representation of $G$ in the case $G= \SL_2$, i.e. the algebra of functions on quantum group $\SL_2$ is generalized by the matrix coefficients $t_{11}, t_{12}, t_{21}, t_{22}$ with the relations 
$$ \begin{array}{ll} t_{11}t_{12} = v^{-2} t_{12}t_{11}, & t_{11}t_{21} = v^{-2} t_{21}t_{11} \\
t_{12}t_{22} = v^{-2} t_{22}t_{12}, & t_{21}t_{22} = v^{-2} t_{22}t_{21} \\
t_{12}t_{21} = t_{21}t_{12} , & t_{11}t_{22} -  t_{22}t_{11} = (v^{-2} - v^2) t_{12}t_{21} \\
	& t_{11}t_{12} - v^{-2}t_{12}t_{21} =1 \end{array}$$
From this presentation of the algebra, L. Korogodski and Y. Soibelman \cite{korosoibel} obtained the description of all the irreducible (infinite-dimensional) unitarizable representations of the quantized algebra of functions $\mathcal F_v(G)$: 
For the particular case of $\mathcal F_v(\SL_2(\mathbb C))$, its complete list of irreducible unitarizable representations consists of:
\begin{itemize}
\item One dimensional representations $\tau_t, t \in \mathbb S^1 \subset \mathbb C$, defined by
$\tau_t(t_{11}) = t, \tau_t(t_{22}) = t^{-1}, \tau_t(t_{12}) = 0, \tau(t_{21}) = 0.$ 
\item  Infinite-dimensional unitarizable $\mathcal F_v(\SL_2(\mathbb C))$-modules $\pi_t, t \in \mathbb S^1$  in $\ell^2(\mathbb N)$,  with an orthogonal basis $\{ e_k\}_{k=0}^\infty$, defined by 
$$\begin{array}{ll}
\pi_t(t_{11}) : & {\left\{ \begin{array}{ll}  e_0 \mapsto 0, &     \\ 
												e_k \mapsto \sqrt{1-v} e_{k-1}, & k > 0, \end{array}\right. }\\
\pi_t(t_{22}) : & {\begin{array}{ll} e_k \mapsto \sqrt{1-v} e_{k+1}, & k \geq 0 \end{array},}\\
\pi_t(t_{12}) : & {\begin{array}{ll} e_k \mapsto tv^{2k} e_{k}, & k \geq 0 \end{array},}\\
\pi_t(t_{21}) : & {\begin{array}{ll} e_k \mapsto t^{-1}v^{2k+1} e_{k}, & k \geq 0 \end{array}.}
\end{array}$$
\end{itemize}

For the general case of $\mathcal F_v(G)$ , consider the algebra homomorphism $\mathcal F_v(G) \to \mathcal F_v(\SL_2(\mathbb C))$, dual to the canonical inclusion $\SL_2(\mathbb C) \hookrightarrow G_{\mathbb C}$.  Then every irreducible unitarizable representation of the quantized algebra of functions $\mathcal F_v(G)$ is equivalent to one of the representations from the list:
\begin{itemize}
\item The representations $\tau_t, t = \exp(2\pi \sqrt{-1} x)\in \mathbb T= \mathbb S^1$, $$\tau_t(C^\lambda_{\nu, s; \mu, r}) = \delta_{r,s} \delta_{\mu, \nu}  \exp(2\pi \sqrt{-1} \mu(x)). $$
\item The representations $\pi_i = \pi_{s_{i_1}} \otimes \dots \otimes \pi_{s_{i_k}},$ if $w= s_{i_1}... s_{i_k}$ is the reduced decomposition of the element $w$ into a product of reflections, where the representations $\pi_{s_i} $ is the composition of the homomorphisms 
$$ \mathcal \pi_{s_i}=\pi_{-1}\circ p : F_v(G) \twoheadrightarrow \mathcal F_v(\SL_2(\mathbb C)) \longrightarrow \End
\ell^2(\mathbb N)$$ 
\end{itemize}
The purpose of this paper is to obtain the same kind results for the quantized algebras of functions on affine Hecke algebras and quantum Schur-Weyl algebras. 
[Remark that it should be more reasonable to name them as the quantized algebras of functions on {\it quantum affine Weyl groups}, but the ``non-affine counterpart" - {\it the quantum Weyl group} terminology was reserved by L. Korogodski and Y. Soibelman for some objects of different kind - the algebras generated not only by the quantized reflections but also the quantized universal algebra.]
\par
We start from the following fundamental remarks:
\begin{itemize}
\item The affine Hecke algebras $\mathbb{H}(v,W^r_{aff})$ and the $v$-Schur algebras $\mathbb S_{n,r}(v)$ are in a complete Schur-Weyl duality. It is therefore easy  to conclude that the corresponding quantized algebras of functions, what we are going to define are also in a complete Schur-Weyl duality.
\item The negative universal enveloping algebras $U^{-}_v(\hat{\mathfrak{sl}}_n) \otimes_A R $, where $A = \mathbb C[v,v^{-1}]$, $R$ is the center of  $U^{-}_v(\hat{\mathfrak{sl}}_n)$, is isomorphic to the Hall algebras $U^{-}_n(\mathbb{H}(v, W^r_{\aff}))$ and there is a natural map $\Theta$ from the last onto the $v$-Schur algebra $\mathbb S_{n,r}(v)$. From this we then have some maps between the quantized algebras of functions
$$\mathcal F(\mathbb S_{n,r}(v)) \to \mathcal F(U^{-}_v(\hat{\mathfrak{sl}}_n)) \to \mathcal F(U_v(\hat{\mathfrak{sl}}_n)).$$ The irreducible representations of $\mathcal F(\mathbb{S}_{n,r}(v))$ could be found in the set of restrictions of irreducible untarizable representations of the quantized algebras $\mathbb C[\SL_r]_q, 0<q<1$, of functions on the quantum group of type ${\SL}_r$.
\item For complex algebraic groups $G$ the irreducible unitarizable $\mathbb C[G]_q$-modules are completely described \cite{korosoibel} for $0<q<1$.
\end{itemize}
\par
Our main result describes the complete set of irreducible unitarizable $\mathcal F(\mathbb{H}(v,W^r_{\aff}))$-modules and $\mathcal F_v(\mathbb S(n,d))$-modules, Theorems \ref{thm2.3}, \ref{thm2.4}, \ref{thm3.1}, \ref{thm3.2}.
\par
Let us describe the paper in more detail: Section 1 is a short introduction to the related subjects and we define the quantized algebras of functions $\mathcal F_v(W^r_{\aff}) = \mathcal F(\mathbb H(v,W^r_{\aff}))$ and $\mathcal F_v(\mathbb S(n,d)):= \mathcal F(\mathbb S_{n,r}(v))$. In \S2 we give a full description of all irreducible unitarizable representations of $\mathcal F_v(W^r_{\aff})$. In \S3 we do the same for the $v$-Schur algebras $\mathcal F(\mathbb S_{n,r}(v))$.
\par
{\sc Notation.}
Let us fix some conventions of notation. Denote $F$ a ground local field of characteristic $p$, $\mathbb C$ the field of complex numbers, $\mathbb Z$ the ring of integers, $\SL_r$ the special linear groups of matrices of sizes $r \times r$ with determinant $1$, $\mathbb G$ an algebraic group, $G=\mathbb G(F)$ the group of rational $F$-points, $\mathbb T$ some maximal torus in $G$, $X^*(\mathbb T)$ the root lattice, $X_*(\mathbb T)$ the co-root lattice, $\mathbb C[G]_q$ the quantized algebra of complex-valued functions on quantum group associated to $G$, $\mathbb S_{n,r}(q)$ the $q$-Schur algebra, $\mathbb S_{n,r}(v)$ the $v$-Schur algebra, $\mathcal F(\mathbb H(v,W^r_{\aff}))$ the quantized algebra of functions on affine Hecke algebra, $\mathcal F(\mathbb S_{n,r}(v))$ the quantized algebra of functions on quantum $v$-Schur algebra. 

\section{Definition of the quantized algebras of functions}
We introduce in this section the main objects - the quantized algebras $\mathcal F_v(W^r_{\aff})$ of functions on affine Hecke algebras. As remarked in the introduction, it should be better to name the quantized algebras  of functions on quantum affine Weyl groups, but we prefer in this paper this terminology in order to avoid any confusion with the terminology from L. Korogodski and Y. Soibelman \cite{korosoibel}.

\subsection{} 
\subsubsection{p-adic presentation}
Let us first recall the definition of Iwahori-Hecke algebras. Let $F$ be a $p$-adic field, i.e. a finite extension of $\mathbb Q_q$, which is by definition the completion with respect to the ultra-metric norm of the rational field of the ring $\mathbb Z_p := \varprojlim \mathbb Z/p^n\mathbb Z$. Denote $\mathcal O$ the ring of integers in $F$, $\mathcal O^\times$ the group of units in $\mathcal O$, $G=\SL_2(F)$, $B = \{ \left[ \begin{array}{cc} x & y \\ 0 & x^{-1} \end{array}\right] ; x,y \in F , x \ne 0 \}$ the Borel subgroup of $G$, $T = \{ \left[ \begin{array}{cc} x & 0 \\ 0 & x^{-1} \end{array}\right] ; x\in \mathcal O^\times \}$ the maximal torus, and $N= \{ \left[ \begin{array}{cc} 1 & y \\ 0 & 1 \end{array}\right] ; y \in \mathcal O \}$ the unipotent radical of $G$. It is easy to check that $B=TN$. Define the so called {\sl Iwahori-Hecke subgroup}
$$I = \{ \left[ \begin{array}{cc} x & y \\ \varpi z & w \end{array}\right] ; x,w \in \mathcal O^\times , y,z \in \mathcal O \},$$  where $\varpi$ is the generic presentative in the presentation of the principal ideal $\mathcal P = \varpi \mathcal O$.
Let us denote $\mu(x)$ the Haar measure on the locally compact group $G=\SL_2(F)$, $\mu(I)= vol(I)$ the volume of $I$ with respect to this Haar measure, $\chi_I$ the characteristic function of the set $I$, $e_I := \frac{1}{\mu(I)} \chi_I$ the idempotent, $e_I^2 = e_I = e_I$, defining a projector. The {\sl Iwahori-Hecke algebra} $\mathbf{IH}(G,I)$ is defined as $$\mathbf{IH}(G,I) = e_I\mathbf H(G)e_I = $$ $$=\{ f: \SL_2(F) \to \mathbb C ; f(hxk) \equiv f(x), \forall h,k\in I, f\in \mathbf H(G) := C^\infty_c(G) \},$$ where $\mathbb H(G) := C^\infty_c(G) $ is the involutive algebra of smooth (i.e. locally constant) functions on $\SL_2(F)$ with compact support, with the well-known convolution product 
$$(f * g)(x) := \int_G f(y) g(y^{-1}x) d\mu(y)$$ and involution as usually.
Recall that the affine Weyl group $W^1_{\aff}$ is defined as $\tilde{W} /\{\pm 1\}$, where $\tilde{W} = \langle D, Dw \rangle$, the group generated by two generators $w := \left[ \begin{array}{cc} 0 & 1 \\ -1 & 0 \end{array} \right]$ and $D:= \left[ \begin{array}{cc} \varpi & 0 \\ 0 & \varpi^{-1} \end{array}\right]$. It is coincided with the dihedral group. Let us choose the following generators $w_1 = w = \left[\begin{array}{cc} 0 & 1 \\ -1 & 0\end{array}\right]$ and $w_2 = \Pi := \left[ \begin{array}{cc} 0 & -\varpi^{-1} \\ \varpi & 0 \end{array}\right]$.  It is well-known the relations $$ w_1w_2w_1^{-1} = w_2^{-1}, \quad w_1^2 = -1,\quad w_2^2 = -1 ,$$ or the standard braid relations
$$ w_1w_2w_1 = w_2w_1w_2, \quad w_1^2 = -1,\quad w_2^2 = -1 .$$ The group $W^1_{\aff}$ is discrete and infinite, and every element of $W^1_{\mathrm aff}$ can be presented as a reduced word in $w_1$ and $w_2$ , namely $w = w_{i_1}\dots w_{i_k}$. The group $G$ can be presented as the union of the double coset classes $G = I.W^1_{\aff}.I$. Let us denote $f_w$ the characteristic function of the coset class $IwI, w\in W_{\aff}$. If $w = w_{i_1} \dots w_{i_k}$ is a reduced presentation of $w\in W^1_{\aff}$ then $f_{w_{i_1}} * \dots * f_{w_{i_k}}$ is independent of the reduced presentation of $w$ and $f_w =f_{w_{i_1}} * \dots * f_{w_{i_k}}$. Let us denote $f_i = f_{w_i}, i=1,2$. We have therefore a correspondence
$$ w\in W^1_{\aff} \mapsto f_w\in \mathbf {IH}(G,I),$$ subject to the relations
$$\left\{ \begin{array}{rcl} f_if_jf_i & = & f_jf_if_j \\
f_i^2 & = & (q-1)f_i + q, \mbox{ with } q=(\mathcal O:\mathcal P).\end{array}\right.$$ 
Let us do a change of variable $v := \frac{1}{\sqrt{q}}$ then we have
$$\left\{ \begin{array}{rcl} f_if_jf_i & = & f_jf_if_j \\
(f_i+1)(f_i-v^{-2}) &= & 0.\end{array}\right.$$ 
This is the so called {\sl Coxeter presentation} of the Iwahori-Hecke algebra in $\SL_2$ case. 
\par
For rank $r$ groups of type $A$, i.e. $\SL_r$ we have the same picture, see for example \cite{chrisginzburg}. Let us consider  also the Hecke algebra $\mathbf H(G) = C^\infty_c(G)$, of smooth (i.e. locally constant) functions on $G$ with compact support, under convolution product and involution. Corresponding to the map of rings $$\mathbb F_q \longleftarrow \mathcal O \longrightarrow F,\mbox{ with } q=p^\ell = (\mathcal O : \mathcal P), \mbox{ for some integer } \ell$$ we have the maps of the groups of rational points $$\mathbb G(\mathbb F_q) \longleftarrow \mathbb G(\mathcal O) \longrightarrow \mathbb G(F).$$ 
The preimage in $\mathbb G(\mathcal O)$ of the Borel subgroup $B(\mathbb F_q)$ is called the {\sl Iwahori subgroup}. It was shown that $G= \mathbb G(F) = I.W^r_{\aff} .I$. The Iwahori-Hecke algebra $\mathbf{IH}(G,I)$ is defined as the algebra of smooth $I$-bi-invariant functions with compact support on $\mathbb G(F)$ under convolution and involution as a sub-algebra of the Hecke algebra $\mathbf H(G)= C^\infty_c(G)$. Denote by $f_{s_i}$ the characteristic function of the double coset class $I.s_i.I$ in $G= \cup_{w\in W^r_{\aff}} I.w.I$, and normalize as in the rank one case we also obtain the relations
$$f_{s_i} f_{s_j} f_{s_i} = f_{s_j} f_{s_i} f_{s_j}, $$
$$(f_{s_i} +1)(f_{s_i} - v^{-2}) = 0, \mbox{ where } v = \frac{1}{\sqrt{q}}, q = p^l = (\mathcal O : \mathcal P),$$
$$f_\sigma f_\gamma = f_{\sigma \gamma} \mbox{ if } \ell(\sigma\gamma) = \ell(\sigma) + \ell(\gamma).$$

\subsubsection{Affine Hecke algebras $\mathbb H(v,W^r_{\aff})$} As usually let us denote $v$ the formal quantum parameter.
{\it (Abstract) Iwahori-Hecke algebras or affine Hecke} are defined in two equivalent ways: in {\sl Coxeter presentation} as group algebras of affine Weyl groups and in {\sl Bernstein presentation } as some abstract algebras presented by generators with relations. 
In Coxeter presentation:
\begin{definition}
An (abstract) Iwahori-Hecke or affine Hecke algebra {\rm is an $\mathbb C[v,v^{-1}]$-algebra generated by $T_\sigma, \sigma \in W^r_{\aff}$, subject to the relations:}
$$T_{s_i} T_{s_j} T_{s_i} = T_{s_j} T_{s_i} T_{s_j}, $$
$$(T_{s_i} +1)(T_{s_i} - v^{-2}) = 0,$$
$$T_\sigma T_\gamma = T_{\sigma \gamma} \mbox{ if } \ell(\sigma\gamma) = \ell(\sigma) + \ell(\gamma).$$
\end{definition}

Let us go to the Bernstein presentation of affine Hecke algebras as some abstract algebras presented by generators with relations. 
\begin{definition}
An affine Hecke algebra in Bernstein presentation {\rm is an $\mathbb C[v,v^{-1}]$-algebra with generators $T^{\pm}_i, i=1,\dots r-1, $ and $X^{\pm}_j, j=1,\dots, r$, subject to the relations:}
$$\begin{array}{ll}
T_iT_i^{-} = 1 = T^{-}_i T_i,\quad & (T_i + 1)(T_i - v^{-2}) = 0,\\
T_iT_{i+1}T_i = T_{i+1}T_iY_{i+1},\quad & T_iT_j = T_jT_i, \mbox{ if } \vert i-j\vert > 1,\\
X_iX_i^{-} = 1 = X_i^{-} X_i, \quad & X_iX_j = X_j X_i,\\
T_iX_iT_i = v^{-2} X_{i+1}, \quad & X_jT_i = T_iX_j, \mbox{ if } J \ne i,i+1.
\end{array}$$
\end{definition}
In this definition we denoted $T_i$ in place of $T^{+}_i$, $X_i$ in place of $X^{+}_i$, etc.... we keep this agreements in the future use.

The isomorphism between two definitions can be established as follows. Associate $T_{s_i} \mapsto T_i$ and $\tilde{T}^{-1}_\sigma \mapsto X_1^{\sigma_1}\dots X_r^{\sigma_r}$, where  $\tilde{T}_\sigma := v^{\ell(\sigma)} T_\sigma$, if $\sigma = (\sigma_1,\dots,\sigma_r)$ is dominant.

\subsection{}
\subsubsection{Admissible representations of $p$-adic groups}
Let us recall that a representation of $p$-adic group is called {\it supercuspidal} iff all its matrix coefficients have compact support modulo the center of the group. It is well-known the following fact: 
Given any irreducible representation $\pi$ of $G$, there exists a Levi subgroup $L$ and a supercuspidal representation $\sigma$ of $L$ such that $\pi$ is a sub-quotient of the induced representation $\imath^G_P(\sigma) := \ind_P^G \infl \sigma$.
Every representation of the form $\imath^G_P(\sigma)$ has finite length for any irreducible representation of $P$ and the other pair $(L',\sigma')$ has the same properties as $(L,\sigma)$ if and only if there exists an element $x\in G$ such that $L'= xLx^{-1}$ and $\sigma' = \sigma^x, $ where $\sigma^x(h) := \sigma(xhx^{-1})$. The pair $(L,\sigma)$ is called a {\it cuspidal pair} and the conjugacy class of $(L,\sigma)$ is called the {\it support} of $\pi$. Two pairs $(L,\sigma)$ and $(L',\sigma')$ are called {\it innertially equivalent} iff there exist $x\in G$ and $\chi\in \mathcal X^{unr}_L$ such that $L'=xLx^{-1}$ and $\sigma' = (\sigma \otimes \chi)^x$. Given an innertially equivalent class $s = \langle (L,\sigma)\rangle$ one defines the sub-category $\mathcal R^s(G)$ of the category $\mathcal R(G)$ of smooth representations, consisting  of all representations, all the sub-quotients of which have support in $s$. One of the well-known result of Bernstein is the fact that $\mathcal R(G) = \times_s \mathcal R^s(G)$ as the direct product of categories. The category $\mathcal R_{cusp}(G) := \times_s \mathcal R^s(G)$. Another well-known result of I. Bernstein, A. Borel and P. Kutzko is the fact that there is an equivalence from the category of unramified representations $\mathcal R^{unr}(G)$, for  $G=\SL_2$,  to the category of finite dimensional representations of the Iwahori-Hecke algebra $\mathbb H(G,I)$. The general case was treated in numerous  works, see for example,  Henniart \cite{henniart}.

\subsubsection{Dipper-James construction of irreducible finite dimensional representations of $\mathbb H(v,W^r_{\aff})$} 
For affine Hecke algebras of type $A_{r-1}$ there are constructions of all irreducible finite dimensional representations parametrized by Young tableaux, or partitions. Let us recall it in brief form. For each Young diagram $\lambda$ a so called Specht $\mathbb H(v,W^r_{\aff})$-module $S^\lambda$ was defined in \cite{dipperjames} and for the value $v=q$ not a root of unity provide a complete list of irreducible finite dimensional representations of $\mathbb H(q,W^r_{\aff})$ modules.

If $q$ is a primitive $\ell$th root of unity, Dipper and James \cite{dipperjames} constructed also a complete set of $\mathbb H(q,W^r_{\aff})$ modules $D^\mu$, parametrized though all Young diagram with at most $\ell-1$ rows of equal length.  Let us describe this construction in more detail.  Let $\lambda = (\lambda_1, \dots, \lambda_r)$, $Y_\lambda = \mathfrak S_{\lambda_1} \times \dots \times \mathfrak S_{\lambda_r} \subset \mathfrak S_n$. Define the symmetrization 
$$\Sym_\lambda := \sum_{w \in Y_\lambda} T_w$$ and the anti-symmetrization
$$A_\lambda := \sum_{w \in Y_\lambda} (-q)^{n(n-1)/2-\ell(w)}T_w.$$
Let $S^\lambda$ be the submodule of the induced $\mathbb H(q, W^r_{\aff})$ module $W^\lambda \cong \mathbb H(q, W^r_{\aff}) \otimes_{\mathbb H(\lambda)} \mathbb C$,
[where $\mathbb H(\lambda)$ is the sub-algebra generated by $T_i$ such that $s_i \in Y_\lambda$], generated by $A_\lambda' W^\lambda$ for $\lambda'$ is obtained from $\lambda$ by interchanging rows with columns, $$S^\lambda = \mathbb H(q, W^r_{\aff}) A_{\lambda' }\mathbb H(q,W^r_{\aff}) \Sym_\lambda.$$ It was proven that there exists an explicit basis of the $\mathbb H(v,W^r_{\aff})$ modules
$$\mathbb H(v,W^r_{\aff}) A_{\lambda'}\mathbb H(v,W^r_{\aff}) \Sym_\lambda \subseteq \mathbb H(v,W^r_{\aff}),$$ which is evaluable at $q \in \mathbb C^\times$ and such that the basis elements evaluated at $q$ remain linearly independent over $\mathbb C$ for all $q\in\mathbb C^\times$. 
Let $(.,.)$ be the bilinear form on the $\mathbb H(v,W^r_{\aff})$ module $W^\lambda$. The the modules $D^\lambda = S^\lambda /(S^\lambda \cap (S^\lambda)^\perp)$ are either 0 or simple. The Young diagram is called $\ell$-regular iff it has at most $\ell-1$ rows of equal length. The module $D^\mu$ is nonzero if and only if $\mu$ is $\ell$-regular.
We refer the reader to the original work of Dipper and James \cite{dipperjames} for a detailed exposition.

\subsubsection{The Langlands correspondence} Recall that a representation of $p$-adic group $G$ is called smooth if the stabilizer of any vector is an open-closed subgroup in $G$. Let us denote $\tilde{V}$ the contragradient representation of $V$, Let $\rho: G=\mathbb G(F)\to \End V$ be an admissible (i.e. smooth and $\tilde{\tilde{V}} = V$) representation of $G$. One of the most important properties of admissible representations of $p$-adic groups is the fact that the space $V^I$ of $I$-invariant vectors in an admissible representation $V$, is finite dimensional. For every element $f$ from the Iwahori-Hecke algebra $\mathbf{IH}(G,I)\cong C^\infty_c(I\setminus G / I)$ we associate an operator in finite dimensional vector space $V^I$,
$$\rho(f) := \int_{\mathbb G(F)} f(x) \rho(x) dx .$$ It is not hard to see that this correspondence gives us a representation of the Iwahori-Hecke algebra $\mathbf{IH}(G,I)$ in the finite dimensional space $V^I$. It was proven that the correspondence $V \mapsto V^I$ provides a functor from, and is indeed an equivalence between the category of admissible representations of $G$ generated by $I$-fixed vectors and the category of finite dimensional representations of the Iwahori-Hecke algebra $\mathbf{IH}(G,I) \cong \mathbb H(v,W^r_{\aff})\vert_{v = q}$. 
This result was essential proven by A. Borel, P. Kutzko end Bernstein in rank one case and by Harris-Taylor \cite{harristaylor} and Henniart \cite{henniart} in the general (rank $r$) case. We refer the readers to  \cite{harristaylor} and \cite{henniart} for more detailed exposition of the {\sl local Langlands Correspondence}.

\subsection{}
We can define now our main objects - the quantized algebras of functions on quantum affine Hecke algebras.
\subsubsection{Quantized algebras of functions}
Let us consider the product of matrix coefficients, associated with the product of elements of the affine Hecke algebra,  of finite dimensional representations, see \cite{lusztig}. With respect to this product we have some non-commutative algebras.

\begin{definition}{\rm The} quantized algebra $\mathcal F(\mathbb H(v,W^r_{\aff}))$ or $\mathcal F_v(W^r_{\aff})$ of functions on the quantum affine Hecke algebra  $\mathbb H(v,W^r_{\aff})$ {\rm is by definition the algebra generated by matrix coefficients of all finite-dimensional representations of the quantum affine Hecke algebra $\mathbb H(v,W^r_{\aff})$
}\end{definition}
\subsubsection{Inclusion}
\begin{proposition} The natural inclusion $W^1_{\aff} \hookrightarrow W^r_{\aff}$ induces a natural projection of quantized algebras of functions
$$\mathcal F(\mathbb H(v,W^r_{\aff}))\twoheadrightarrow \mathcal F(\mathbb H(v,W^1_{\aff}))  .$$ 
\end{proposition}
{\sc Proof.} It easy an easy consequence from the corresponding inclusion of the affine Weyl groups, $W^1_{\aff} \hookrightarrow  W^r_{\aff}$ . \hfill $\square$

\section{Irreducible representations}
The main subject of this section is to describe all (up to unitary equivalence) inequivalent unitarizable  representations of the quantized algebras of functions on affine Hecke algebras. We describe first in the rank 1 case and then use the projection $\mathcal F(\mathbb H(v,W^r_{\aff}))\twoheadrightarrow \mathcal F(\mathbb H(v,W^1_{\aff})) $ to maintain the general case.

\subsection{Rank 1 case}
\begin{lemma}\label{lemma2.1}
The quantized algebra ${\mathcal F}_v(W^f \ltimes X_*(T))$ is generated by the restrictions $t_{11}\vert_{W^f \ltimes X_*(T))}$ and $t_{12}\vert_{W^f \ltimes X_*(T))}$ with some defining relations.
\end{lemma}
{\sc Proof.}
It was proven in L. Korogodski and Y. Soibelman \cite{korosoibel} that in every finite-dimensional representation of $\mathcal F[\SL_2(\mathbb C)]_q$, there exists an action of quantum Weyl elements $\bar{w}$. For the groups of type $A_1$, the root and coroot lattices are isomorphic $X^*(\mathbb T) \cong X_*(\mathbb T)$. We can therefore see $W_{\aff} = W^f \ltimes X^*(\mathbb T)\cong \bar{W}_{\aff} = W^f \ltimes X_*(\mathbb T)$ as some subgroups of $\SL_2(\mathbb C)$. Therefore we have the restrictions of the representations from the list of irreducible representations of $\SL_2(\mathbb C)$. Two generators of $W^1_{\aff}$ are $w = \left[ \begin{array}{cc} 0 & 1\\ -1 & 0\end{array} \right]$ and $D= \left[\begin{array}{cc} \varpi & 0\\ 0 & \varpi \end{array} \right]$. In the representation described in \cite{korosoibel}, they are defined by two matrix elements $t_{11}$ and $t_{12}$, restricted to our affine Weyl group.
\hfill $\square$
\begin{lemma} \label{lemma2.2}
Every irreducible unitarizable representation of ${\mathcal F}_v(W^f \ltimes X_*(T))$  can be obtained by restricting some irreducible unitarizable representations of ${\mathcal F}_v(\SL_2(\mathbb C)$.
\end{lemma}
{\sc Proof.} First remark that if $V$ is a representation of ${\mathcal F}_v(W^f \ltimes X_*(T))$ and $\ind V={\mathcal F}_v(\SL_2(\mathbb C)) \otimes_{{\mathcal F}_v(W^f \ltimes X_*(T)) } V $ the induced representation of ${\mathcal F}_v(\SL_2(\mathbb C))$, then there is the well-known Frobenius duality
$$\Hom(\ind V, W) \cong \Hom(V,W\vert_{{\mathcal F}_v(W^f \ltimes X_*(T))}) .$$
Let us consider a ${\mathcal F}_v(W^f \ltimes X_*(T))$  module $V$.
Taking induction $\ind V = {\mathcal F}_v(\SL_2(\mathbb C)) \otimes_{{\mathcal F}_v(W^f \ltimes X_*(T)) } V$ , we have a ${\mathcal F}_v(\SL_2(\mathbb C)) $ module. The irreducible ones can be therefore obtained from the list of irreducible unitarizable reprenatations $\pi_{w,t}$ of ${\mathcal F}_v(\SL_2(\mathbb C))$. 
\hfill $\square$

Let us denote the restrictions of representations of $\mathcal F[\SL_2]_q$ on $\mathcal F(\mathbb H(v,W^r_{\aff}))$ by the same letters.
\begin{theorem} \label{thm2.3} Every irreducible unitarizable representations of $\mathcal F_v(W^1_{\aff})$ is equivalent to one of the unitarily inequivalent representation from the list:
\begin{enumerate}
\item The representations $\tau_t, t\in \mathbb S^1$, defined by 
$$\begin{array}{cc}  \tau(t_{11}) = t, & \tau(t_{22}) = t^{-1},\\
 								 \tau(t_{21}) = 0, & \tau(t_{12}) = 0,\end{array}$$

\item The representations $\pi_{w,t}, w\in W^f, t\in \mathbb S^1$, defined by 
$$\begin{array}{ll}
\pi_t(t_{11}) : & {\left\{ \begin{array}{ll}  e_0 \mapsto 0, &     \\ 
												e_k \mapsto \sqrt{1-v} e_{k-1}, & k > 0, \end{array}\right. }\\
\pi_t(t_{22}) : & {\begin{array}{ll} e_k \mapsto \sqrt{1-v} e_{k+1}, & k \geq 0 \end{array},}\\
\pi_t(t_{12}) : & {\begin{array}{ll} e_k \mapsto tv^{2k} e_{k}, & k \geq 0 \end{array},}\\
\pi_t(t_{21}) : & {\begin{array}{ll} e_k \mapsto t^{-1}v^{2k+1} e_{k}, & k \geq 0 \end{array}.}
\end{array}$$
\end{enumerate}
\end{theorem}
{\sc Proof.} It is directly deduced from Lemmas \ref{lemma2.1}, \ref{lemma2.2} and the following fact.
Let us now recall  that L. Korogodski and Y. Soibelman \cite{korosoibel} obtained the description of all the irreducible (infinite-dimensional) unitarizable representations of the quantized algebra of functions $\mathcal F_v(G)$: 
For the particular case of $\mathcal F_v(\SL_2(\mathbb C))$, its complete list of irreducible unitarizable representations consists of:
\begin{itemize}
\item One dimensional representations $\tau_t, t \in \mathbb S^1 \subset \mathbb C$, defined by
$\tau_t(t_{11}) = t, \tau_t(t_{22}) = t^{-1}, \tau_t(t_{12}) = 0, \tau(t_{21}) = 0.$ 
\item  Infinite-dimensional unitarizable $\mathcal F_v(\SL_2(\mathbb C))$-modules $\pi_t, t \in \mathbb S^1$  in $\ell^2(\mathbb N)$,  with an orthogonal basis $\{ e_k\}_{k=0}^\infty$, defined by 
$$\begin{array}{ll}
\pi_t(t_{11}) : & {\left\{ \begin{array}{ll}  e_0 \mapsto 0, &     \\ 
												e_k \mapsto \sqrt{1-v} e_{k-1}, & k > 0, \end{array}\right. }\\
\pi_t(t_{22}) : & {\begin{array}{ll} e_k \mapsto \sqrt{1-v} e_{k+1}, & k \geq 0 \end{array},}\\
\pi_t(t_{12}) : & {\begin{array}{ll} e_k \mapsto tv^{2k} e_{k}, & k \geq 0 \end{array},}\\
\pi_t(t_{21}) : & {\begin{array}{ll} e_k \mapsto t^{-1}v^{2k+1} e_{k}, & k \geq 0 \end{array}.}
\end{array}$$
\end{itemize}

\hfill $\square$

\subsection{Rank r case}
Let us consider the representations which are the 
composition of the homomorphisms 
$$ \mathcal \pi_{s_i}= \pi_{-1} \circ p : F_v(G) \twoheadrightarrow \mathcal F_v(\SL_2(\mathbb C)) \longrightarrow \End 
\ell^2(\mathbb N), $$ 

\begin{theorem} \label{thm2.4}
Every irreducible unitarizable representation of $\mathcal F_v(W^r_{\aff})$ is equivalent to one of the unitarily inequivalent representations:
\item The representations $\pi_{w,t}= \pi_{s_{i_1}} \otimes \dots \pi_{s_{i_1}} \otimes \tau_t$, $w= s_{i_1} \dots s_{i_k} \in W^f$ is a reduced decomposition of $w$, $ t\in \mathbb S^1$. 
\end{theorem}
{\sc Proof}. 
For the general case of $\mathcal F_v(G)$ , consider the algebra homomorphism $\mathcal F_v(G) \twoheadrightarrow \mathcal F_v(\SL_2(\mathbb C))$, dual to the canonical inclusion $\SL_2(\mathbb C) \hookrightarrow G_{\mathbb C}$.  Then every irreducible unitarizable representation of the quantized algebra of functions $\mathcal F_v(G)$ is equivalent to one of the representations from the list:
\begin{itemize}
\item The representations $\tau_t, t = \exp(2\pi \sqrt{-1} x)\in \mathbb T= \mathbb S^1$, $$\tau_t(C^\lambda_{\nu, s; \mu, r}) = \delta_{r,s} \delta_{\mu, \nu}  \exp(2\pi \sqrt{-1} \mu(x)).$$ 
\item The representations $\pi_i = \pi_{s_{i_1}} \otimes \dots \otimes \pi_{s_{i_k}},$ if $w= s_{i_1}... s_{i_k}$ is the reduced decomposition of the element $w$ into a product of reflections, where the representations $\pi_{s_i} $ is the composition of the homomorphisms 
$$ \pi_{s_i}= \pi_{-1} \circ p :  \mathcal  F_v(G) \twoheadrightarrow \mathcal F_v(\SL_2(\mathbb C))  \longrightarrow \End 
\ell^2(\mathbb N), $$ 
\end{itemize}

\hfill$\square$

\section{Schur-Weyl duality}
The Schur-Weyl duality is well-known for finite-dimensional representations of quantum affine Hecke algebras and quantum $v$-Schur algebras. For (possibly infinite dimensional) representations of the quantized algebras of functions on them we also have this kind of duality. We use it then to describe (possibly infinite-dimensional) representations of $q$-Schur algebras. The main idea is to use the maps
$$U_n(\hat{sl}_r) \twoheadrightarrow U^{-}_n(\hat{sl}_r) \rightarrowtail U^{-}_n(\hat{sl}_r) \otimes_A R \twoheadrightarrow \mathbb S_{n,r}(v)$$

\subsection{}
\subsubsection{$v$-Schur algebras $\mathbb S_{n,r}(v)$} 
We recall first the definition of the affine $v$-Schur algebras. Let $s\in \mathbb N$ be an nonnegative integer, and $r\in \mathbb N^* = \mathbb N \setminus \{ 0\}$ a positive integer. Denote 
$$\mathcal A^n_r = \{(i_1,\dots ,i_r) ;  1\leq i_1\leq \dots \leq i_r \leq n \}$$ be the fundamental domain of the both actions of $W^r_{\aff} = \hat{\mathfrak S}_r$ on ${\mathbb Z^r}$ on the left by 
$$s_j .(i_1, \dots, i_r) := (i_1, \dots, i_{j+1}, i_j, \dots, ir), 1\leq j < r,$$
$$\lambda . (i_1, i_r) := (i_1 +s\lambda_1, \dots i_r + s\lambda_r ), \lambda \in \mathbb Z^r$$
and on the right by
$$(i_1, \dots, i_r).s_j  := (i_1, \dots, i_{j+1}, i_j, \dots, ir), 1\leq j < r,$$
$$(i_1, i_r).\lambda   := (i_1 +s\lambda_1, \dots i_r + s\lambda_r ), \lambda \in \mathbb Z^r.$$
For an element $\mathbf i \in \mathcal A^n_r$, denote the stabilizer as $\mathfrak S_{\mathbf i}$. Let us consider the projector $e_{\mathbf i}:= \sum_{\delta \in \mathfrak S_{\mathbf i}} T_\delta$. Define the affine $v$-Schur algebra as
$$\mathbb S_{n,r}(v) := \bigoplus_{\mathbf i, \mathbf j \in \mathcal A^n_r} \mathbb H_{\mathbf i,\mathbf j} = \bigoplus_{\mathbf i, \mathbf j \in \mathcal A^n_r} e_{\mathbf i} \mathbb H(v,W^r_{\aff}) e_{\mathbf j}.$$
It was proven that $\mathbb H_{\mathbf i, \mathbf j} = e_{\mathbf i} \mathbb H(v,W^r_{\aff}) e_{\mathbf j}$ is exactly the $\mathbb C[v,v^{-1}]$-linear  span of the element $T_{\sigma} = \sum_{\sigma \in \mathfrak S_{\mathbf i} \setminus \hat{\mathfrak S} / \mathfrak S_{\mathbf j}} T_\sigma$. It was proven that this affine $v$-Schur algebra $\mathbb S_{n,r}(v)$ is a quotient of the modified quantum group $\dot U_v^{-}(\hat{gl}_n)$.

\subsubsection{$v$-Schur duality}
One defines $$\mathbb T_{n,r} := \bigoplus_{\mathbf i \in \mathcal A^n_r} e_{\mathbf i}\mathbb H(v,W^r_{\aff}).$$
Define $T_\sigma := \sum_{\delta \in \sigma} T_\delta$, for each coset class $\sigma \in \mathfrak S_{\mathbf i}\setminus \hat{\mathfrak S}_r$, then $\{ T_{\sigma} \}$ form a basis of $\mathbb T(n,r)$.
The algebra $\mathbb H(v,W^r_{\aff})$ acts on $\mathbb T(n,r)$ by multiplication on the right and the algebra $\mathbb S_{n,r}(v)$ acts on $\mathbb T(n,r)$ on the left by multiplication
$$e_{\mathbf i}he_{\mathbf j}. e_{\mathbf k}h' := \delta_{\mathbf j, \mathbf k} e_{\mathbf i}he_{\mathbf j}h',\forall h,h' \in e_{\mathbf i}\mathbb H(v,W^r_{\aff}).$$
{\bf The Schur-Weyl duality} for finite dimensional representations  is as follows.
$$\mathbb S_{n,r}(v) \cong \End_{\mathbb H(v,W^r_{\aff})} \mathbb T(n,r),$$
$$\mathbb H(v,W^r_{\aff}) \cong \End_{\mathbb S_{n,r}(v)} \mathbb T(n,r),$$
Remark that a geometric realization of this Schur-Weyl duality is an important subject in the Deligne-Langlands interplay and was highly developed, see e.g. \cite{chrisginzburg}.

\begin{theorem} \label{thm3.1} The unitarizable $\mathcal F(\mathbb H(v,W^r_{\aff}))$-modules and $\mathcal F(\mathbb S_{n,r}(v))$ modules are in a complete Schur-Weyl duality
$$\mathcal F(\mathbb S_{n,r}(v)) \cong \End_{\mathcal F(\mathbb H(v,W^r_{\aff}))} \mathbb T(n,r),$$
$$\mathcal F(\mathbb H(v,W^r_{\aff})) \cong \End_{\mathcal F(\mathbb S_{n,r}(v))} \mathbb T(n,r),$$

\end{theorem}

{\sc Proof.}
It is enough to recall that the quantum algebras of functions are consisting of matrix coefficients of all finite dimensional representations of the affine Hecke algebras and affine $v$-Schur algebras respectively. 
\hfill$\square$

\subsection{}
\subsubsection{Restriction maps} 
Let us first recall \cite{lusztig} the definition of the so called modified universal enveloping algebras $\dot U(\mathfrak g)$. Denote as before $X^*(T)$ the weight lattice, $X_*(T)$ the co-weight lattice. For each $\lambda', \lambda'' \in X^*(T)$ define
$${}_{\lambda'}U_{\lambda''} := U(\mathfrak g) / (\sum_{\mu\in X_*(T)}( K_\mu - v^{\langle\mu,\lambda\rangle})U(\mathfrak g) + U(\mathfrak g)\sum_{\mu\in X_*(T)}( K_\mu - v^{\langle\mu,\lambda\rangle}))$$ and the natural projection $$U(\mathfrak g) \twoheadrightarrow {}_{\lambda'}U_{\lambda''}.$$ By definition the modified universal enveloping algebra $\dot U(\mathfrak g)$ is the direct sum
$$\dot U(\mathfrak g) := \bigoplus _{\lambda'\in X^*(T), \lambda'' \in X_*(T)} {}_{\lambda'}U_{\lambda''} .$$

The $v$-Schur algebras can be considered as some quotient of the modified quantized universal enveloping algebras $\dot {U}_v(\mathfrak g)$ which is different from $U(\mathfrak g)$ replacing $U^0(\mathfrak g)= \mathbb C$ by the direct sum of infinite number of copies of $\mathbb C$, one for each element of the weight lattice $X^*(T)$, see G. Lusztig (\cite{lusztig}, {chapters. 23, 29}).
 It was shown that the category of highest weight finite dimensional representations with weight decomposition of $U(\mathfrak g)$ is equivalent to the category of highest weight representations of $\dot U(\mathfrak g)$, but the algebras $U(\mathfrak g)$ admit also the representations without weight decomposition. 

Recall from the work of Schiffmann. The main idea is related with the maps 
$$U_n(\hat{sl}_r) \twoheadrightarrow U^{-}_n(\hat{sl}_r) \rightarrowtail U^{-}_n(\hat{sl}_r) \otimes_A R \twoheadrightarrow \mathbb S_{n,r}(v)$$
\subsubsection{Description of irreducible representations}
\begin{theorem} \label{thm3.2}
The restrictions of irreducible unitarizable $\mathcal F_v(U_n(\hat{sl}_r))$ modules to $\mathcal F_v(\mathbb S_{n,r}(v))$ give a complete list of irreducible unitarizable $\mathcal F_v(\mathbb S_{n,r}(v))$ modules.
\end{theorem}
{\sc Proof}.
The proof combines Lemmas \ref{lemma2.1}, \ref{lemma2.2} and the following fact. 
 In the particular case of $\mathbb S_{2,d}(v)$ Doty and Giaquinto \cite{dotygiaquinto} have a more presice description: The $v$-Schur-Weyl algebra is just the image of the quantized universal eveloping algebra $U_v(\mathfrak {sl}_2)$ in the $d$-tensor product power of the standard 2-dimensional representation. It is isomorphic to the algebra generated by elements $E, F, K$ and $K^{-1}$ subject to the relations:
\begin{enumerate}
\item[(a)] $KK^{-1} = K^{-1}K = 1,$
\item[(b)] $KEK^{-1} = v^2 E, \quad KFK^{-1} = v^{-2} F,$ 
\item[(c)] $EF-FE = \frac{K - K^{-1}}{v-v^{-1}},$
\item[(d)] $(K-v^d)(K-v^{d-2}) \dots (K-v^{-d+2}) (K-v^{-d}) =0. $
\end{enumerate}
We use again the map $\mathcal F_v(\mathbb S(n,d)) \to \mathcal F_v(\mathbb S(2,d))$ associated with the natural inclusion of the Weyl groups $W^1_{\aff} \hookrightarrow W^r_{\aff}$ 
\hfill$\square$
\begin{remark}{\rm
Denote $$(V_{\varepsilon, t},\Pi_{\varepsilon,t})= \left\{ \begin{array}{ll} (\mathbb C, \tau_t), & \mbox{if } \varepsilon = 0, \\   (\ell_2(\mathbb N), \pi_{w, t}), w\in W^f, &\mbox{if } \varepsilon = 1\end{array} \right.$$ and define $$\widetilde{\mathbb T}_{n,r} := \sum^\oplus_{\varepsilon = 0,1} \int^\oplus_{t\in \mathbb S^1} V^*_{\varepsilon, t} \otimes V_{\varepsilon, t} dt.$$ We have therefore the {\bf Schur-Weyl Duality for unitarizable representations}: Every irreducible unitarizable representation of the quantum affine Hecke algebra $\mathbb H(v,W^r_{\aff})$ is a sub-representation of the representation in the space of $\mathbb S_{n,r}(v)$-invariants $\widetilde{\mathbb T}_{n,r}^{\mathbb S_{n,r}(v)}$ and conversely, every irreducible unitarizable representation of the quantum $v$-Schur algebra $\mathbb S_{n,r}(v)$ is a sub-representation of the representation in the space of $\mathbb H(v,W^r_{\aff})$-invariants $\widetilde{\mathbb T}_{n,r}^{\mathbb H(v,W^r_{\aff}) }$ 
}\end{remark}

\section*{Acknowledgments}
This work was completed during the stay of the author as a visiting
mathematician at the Department of mathematics, The University of Iowa. The author would like to express the deep and sincere thanks to Professor Tuong Ton-That and his spouse, Dr. Thai-Binh Ton-That for their effective helps and kind attention they provided during the stay in Iowa, and also for a discussion about the PBW Theorem and Schur-Weyl duality. The deep thanks are also addressed to the organizers of the Seminar on Mathematical Physics, Seminar on Operator Theory in Iowa and the Iowa-Nebraska Functional Analysis Seminar (INFAS), in particular the professors Raul Curto, Palle Jorgensen, Paul Muhly and Tuong Ton-That for the stimulating scientific atmosphere. The deep thanks are addressed to professors Phil Kutzko and Fred Goodman for the useful discussions
during their seminar lectures on Iwahori-Hecke algebras and their representations.
 
The author would like to thank the University of Iowa for the hospitality and the scientific support, the Alexander von Humboldt Foundation, Germany, for an effective support.

{\noindent{\sc Department of Mathematics, The University of Iowa, 14 McLean Hall, Iowa City, IA 52242-1419, U. S. A.}\\  {\rm Email:}\quad {\tt ndiep@math.uiowa.edu}}\\
and\\
{\noindent{\sc Institute of Mathematics, NCST of Vietnam, P. O. Box 631, Bo Ho 10,000, Hanoi, Vietnam}\\  {\rm Email:}\quad {\tt dndiep@hn.vnn.vn}}
\end{document}